\documentclass[12pt]{amsart}

\makeatletter
\def\LaTeX{\leavevmode L\raise.42ex
    \hbox{\kern-.3em\size{\sf@size}{0pt}\selectfont A}\kern-.15em\TeX}
\makeatother

\newcommand{\BibTeX}{{\rm B\kern-.05em{\sci\kern-.025emb}\kern-.08em\TeX}}

\makeatletter
\def\@currentlabel{2.1}\label{e:dispaa}
\def\@currentlabel{2.21}\label{e:dispau}
\def\@currentlabel{2.22}\label{e:dispav}
\def\@currentlabel{2.23}\label{e:dispaw}
\def\@currentlabel{2.24}\label{e:dispax}
\def\theequation{\thesection.\@arabic\c@equation}
\makeatother

\newcommand{\R}{{\mathbb R}}

\newcommand{\ep}{\epsilon}

\newtheorem{lemma}{Lemma}[section]

\newtheorem{theorem}[lemma]{Theorem}
\newtheorem{corollary}[lemma]{Corollary}

\title[Semilinear Elliptic Equations On A Strip]%
{ On  least Energy Solutions to A Semilinear Elliptic Equation in A Strip}
\author[ Henri Berestycki and Juncheng Wei]{}


 \email{hb@ehess.fr} \email{wei@math.cuhk.edu.hk}

\begin{document}

\maketitle

\centerline{\scshape Henri Berestycki}
\medskip
{\footnotesize \centerline{Ecole des hautes etudes en sciences
sociales} \centerline{CAMS, 54, boulevard Raspail, F - 75006 -
Paris, France}
 }

\medskip

\centerline{\scshape Juncheng Wei }
\medskip
{\footnotesize \centerline{Department of Mathematics, Chinese
University of Hong Kong}\centerline{Shatin, Hong Kong, China}
 }

\bigskip

\begin{center}
\emph{Dedicated to Professor L. Nirenberg on the occasion of his
85$^{th}$ birthday,\\ with deep admiration}

\end{center}

\begin{abstract} We consider the following semilinear elliptic equation on a strip:
\[ \left\{\begin{array}{l}
\Delta u-u + u^p=0 \ \mbox{in} \ \R^{N-1} \times (0, L),\\
u>0, \frac{\partial u}{\partial \nu}=0 \ \mbox{on} \ \partial (\R^{N-1} \times (0, L))
\end{array}
\right.\]
where $ 1< p\leq \frac{N+2}{N-2}$.  When $ 1<p <\frac{N+2}{N-2}$, it is shown that there exists a unique $L_{*} >0$ such that for $L \leq L_{*}$, the least energy solution is trivial, i.e., doesn't depend on $x_N$,  and for $L >L_{*}$, the least energy solution is nontrivial. When $N \geq 4, p=\frac{N+2}{N-2}$, it is shown that there are two numbers $L_{*}<L_{**}$ such that the least energy solution is trivial when $L \leq L_{*}$, the least energy solution is nontrivial when $L \in (L_{*}, L_{**}]$, and the least energy solution does not exist when $L >L_{**}$. A connection with Delaunay surfaces in CMC theory is also made.
\end{abstract}

\section{Introduction}
\setcounter{equation}{0}

In this paper, we consider the following semilinear elliptic equation on a strip
\begin{equation}
\label{1}
\left\{\begin{array}{l}
\Delta u-u +u^p=0 \ \mbox{in} \ \R^{N-1} \times (0, L),\\
  \frac{\partial u}{\partial \nu}=0 \ \mbox{on} \ \partial (\R^{N-1} \times (0, L)), \\
u>0, \ u \in H^1 (\R^{N-1} \times (0, L)).
\end{array}
\right.
\end{equation}
Here we assume that
\[ N\geq 2, 1 <p \leq \frac{N+2}{N-2} \ \mbox{if} \ N\geq 3, \mbox{and} \ 1<p <+\infty \ \mbox{if}\ N=2\]
and $\nu$ is the outer normal derivative.

The motivation of this study stems of the work of Dancer on new solutions to
for the following simple superlinear problem
\begin{equation}
 \label{4}
\Delta u-u+ u^p=0 \ \mbox{in} \ \R^N, u>0, p>1.
\end{equation}

If $u(x) \to 0$ as $ |x| \to +\infty$, then the classical work of Gidas-Ni-Nirenberg \cite{gnn} shows that $u$ must be {\it radially symmetric} with respect to one point and thus (\ref{4}) is reduced to an ODE. On the other hand, there are less known results on solutions to (\ref{4}) which do not decay in all directions. Dancer \cite{d1} first constructed  solutions to (\ref{4}) that are periodic in one direction and decays in all the other directions, via local bifurcation arguments. They form a
one-parameter family of solutions which are periodic in the $z$
variable and originate from the decaying solutions of (\ref{4}) in $\R^{N-1}$. We briefly outline Dancer's idea:  Let $T>0$ be the period and consider
\begin{equation}
 \label{5n}
 \left\{\begin{array}{l}
\Delta u-u+ u^p=0 \ \mbox{in} \ \R^N, u>0, \\
u(x^{'}, x_N+T)= u(x^{'}, x_N), u (x^{'}, x_N) \to 0 \ \mbox{as} \ |x^{'}| \to +\infty
\end{array}
\right.
\end{equation}
where we denote $x^{'}=(x_1,..., x_{N-1})$. Dancer then used $T$ as the bifurcation parameter and found a critical value $T_{1}$ such that for $T=T_{1}$, the linearized problem at the lower dimensional decaying solutions has an eigenvalue zero with eigenfunctions decaying in $x^{'}$. Then using the Crandall-Rabinowitz bifurcation theory, near $T_1$, a new solution (different from lower dimensional solution) bifurcates.

In \cite{DKPW}, these periodic solutions are called {\em Dancer's solutions} and they are the building blocks for more complicated ``$2k$-ends'' solutions. In \cite{mal}, Dancer's solutions are also used to build {\it three ends} solutions to the problem in entire space. In fact, geometrically, Dancer's solutions corresponds to the so-called Delaunay solution in CMC theory \cite{CD}. We will comment on this later.  Therefore it becomes natural  to study the solution structure of (\ref{1}).

Problem (\ref{1}) also arises naturally in the study of some nonlinear elliptic equations in an expanding annuli:
\begin{equation}
\label{2}
\left\{\begin{array}{l}
\Delta u-u +u^p=0 \ \mbox{in} \  B_{R+L} \backslash B_R,\\
 u>0,  \frac{\partial u}{\partial \nu}=0 \ \mbox{on} \ \partial (B_{R+L} \backslash B_R), \\
\end{array}
\right.
\end{equation}
where $R \to +\infty$ and $L$ is fixed. The limiting equation of (\ref{2}) as $R \to +\infty$ becomes (\ref{1}).

We note that the corresponding expanding annuli Dirichlet problem
\begin{equation}
\label{2new}
\left\{\begin{array}{l}
\Delta u-u +u^p=0 \ \mbox{in} \  B_{R+L} \backslash B_R,\\
 u>0,  u=0 \ \mbox{on} \ \partial (B_{R+L} \backslash B_R), \\
\end{array}
\right.
\end{equation}
 has been studied by many authors, see \cite{b}, \cite{cw}, \cite{dy}, \cite{li}, \cite{ls1}, \cite{ls2}, \cite{ms} and the references therein.  We are not aware of any study on (\ref{2}). Note that (\ref{2}) is different  from (\ref{2new}). Indeed,   (\ref{2}) admits solutions that are nonzero and lower dimensional, i.e. don't depend on $x_N$-direction. The issue, therefore, is to understand when and how solutions that are not lower dimensional exist.

By suitable scaling (\ref{1}) becomes
\begin{equation}
\label{5}
\left\{\begin{array}{l}
\Delta u- L^2 u +u^p=0 \ \mbox{in} \ \Sigma:= \R^{N-1} \times (0, 1),\\
  \frac{\partial u}{\partial \nu}=0 \ \mbox{on} \ \partial \Sigma, \ u>0, \ u \in H^1 (\Sigma).
\end{array}
\right.
\end{equation}

Here $L$ is the parameter. In this paper, we consider the existence or non-existence as well as nature of {\em least energy solutions}. More precisely, let
\begin{equation}
\label{cl}
c( L):= \inf_{ u \in H^1 (\Sigma), u \not \equiv 0}  \frac{ \int_\Sigma (|\nabla u|^2 + L^2 u^2)}{ (\int_\Sigma u^{p+1})^{\frac{2}{p+1}}}.
\end{equation}

Our main concerns are:

\noindent
Q1: Is $c(L)$ attained?

\noindent
Q2: Is $c( L)$ attained by a nontrivial solution?

\noindent
Q3: Is the least energy solution nondegenerate?

Here, a trivial solution is understood to mean that the solution
does not  depend on $x_N$. Note that for $ p <\frac{N+1}{N-3}$ if
$N\geq 4$ and $1<p<+\infty$ if $N=2,3$, such trivial solutions of (\ref{5})
always exist. Indeed, a ground state solution
(\cite{BL1}-\cite{BL2}) in $\R^{N-1}$ yields such a ``trivial
solution".

Our main theorem is

\begin{theorem}
\label{1.1}

(1) If $ p < \frac{N+2}{N-2}$ when $N\geq 3$ and $p<+\infty$ when $N=2$, then there exists a unique $L_{*}$ such that for $L \leq L_{*}$, $c(L)$ is attained by a trivial solution and for $L >L_{*}$, $c(L)$ is attained by a nontrivial solution.

(2) There exist $ L_2\geq L_{*}$ such that the least energy solution is unique and nondegenerate for any $L \geq L_2$.

(3) If $N \geq 4, p=\frac{N+2}{N-2}$, then there exists two positive constants $L_{*} <L_{**}$ such that for $L \leq  L_{*}$, $c(L)$ is attained by a trivial solution; for $ L \in (L_{*}, L_{**}]$, $c(L)$ is attained by a nontrivial solution; for $ L>L_{**}$, $c(L)$ is not attained.
\end{theorem}

\noindent
{\bf Remark:} The number $L^{*}$ can be computed as follows: Let $w_0$ be the unique  ground state solution  in $\R^{N-1}$
\begin{equation}
\label{w0}
\Delta w_0- w_0 + w_0^p=0, \ w_0=w_0 (|x^{'}|)>0, w_0 \in H^1 (\R^{N-1}).
\end{equation}
(See \cite{BL1}, \cite{k}.) Let $\lambda_1$ be the unique principal eigenvalue of
\begin{equation}
\Delta \phi-\phi+pw_0^{p-1}\phi=\lambda_1 \phi, \phi \in H^1(\R^{N-1}).
\end{equation}
Then we have
\begin{equation}
L^{*}=\frac{\pi}{\sqrt{\lambda_1}}.
\end{equation}
When $N=2$, we can compute explicitly (see \cite{dkw-cpam})
\begin{equation}
 \lambda_1=\frac{(p-1)(p+3)}{4}.
 \end{equation}
In fact, we can say more about the properties of the minimizers and the asymptotic behaviors of $c(L)$ as $L \to 0$ or $L \to +\infty$.  The asymptotic behavior of the least energy solution when $L \to L_{*}$ is given in the appendix.

 Even though Theorem \ref{1.1} is a purely PDE result, this result has a striking  analogy in the theory of constant mean
curvature (CMC) surface in $\R^3$.

  CMC surfaces in $\R^3$ are  equilibria for the area functional subjected to an
enclosed volume constraint. It arises in many physical  and
variational problems. Over the past two decades  a great deal of progress was achieved in understanding complete
CMC  surfaces and their {\em moduli} spaces. Spheres (zero end)  and round cylinders are the first examples of  CMC surfaces. (See
 Alexandrov's \cite{A}.)  Properly embedded CMC surfaces with
nonzero mean curvature  were classified by Delaunay \cite{CD}. These are CMC  rotation
surfaces, called {\em unduloids} (having genus zero and two ends).
These surfaces are derived from two 1-parameter families: one of the family being unduloids with neck radius
$\tau\in (0, \frac{1}{2}]$ and the other being a family of non-embedded surface called {\em nodoids} that can be parameterized by the neck radius  $\tau\in (0, \infty)$.

In very much an analogous way, our solutions in Theorem \ref{1.1} are parameterized by the length $L$. When $L \to +\infty$, these solutions become spikes at the center and correspond to the Delaunay surface that are obtained when $\tau \to 0$. On the other hand, when $L \to L_{*}$ our solution corresponds to Delaunay solution when $\tau \to \frac{1}{2}$. One good way to think of this analogy  is the level sets of $u$. (See \cite{DKPW} for more explanations.) In \cite{DKPW}, del Pino-Kowalczyk-Pacard-Wei used the least energy
 solution near $L^{*}$ and Toda systems  to build more  complicated even-ended solutions of (\ref{4}) in $\R^2$, while in \cite{mal}, Malchiodi used the least energy  solution near $+\infty$ to build $Y-$shaped solutions.

We conjecture that the least energy solution form  a continuous family as $ L$ goes from $L_{*}$ to $+\infty$.

After the paper was completed, we  learned from Prof. M. Esteban that problem (\ref{5}) is also related to the study of Cafferalli-Kohn-Nirenberg  inequality  and it is studied in the work of Dolbeault, Esteban, Loss and Tarantello \cite{det}.

This paper is organized as follows: we prove (1), (2) and (3) of Theorem \ref{1.1} in Sections 2,3 and 4 respectively. In Appendix A, we prove some technical estimates used in Section 4 while in Appendix B we  study the asymptotic behavior of least energy solutions  when $L \to L_{*}$.

\bigskip
{\bf Acknowledgments:} Part of this research was completed while the first author was visiting the University of Chicago. The research of the second author is
partially supported by a General Research Fund from RGC of Hong Kong. The second author thanks CAMS, EHESS for the hospitality he received during his visit in June, 2009.

\section{The subcritical case: Proof of (1) of Theorem 1.1}
\setcounter{equation}{0}

In this section, we study (\ref{5}) for the subcritical case, i.e., $1<p <\frac{N+2}{N-2}$ if $N\geq 3$ and $1<p<+\infty$ when $N=2$. We begin with

\begin{lemma}
\label{l00}
For any $L_0>0$, there exists a constant $C$, independent of $L \leq L_0$, such that for any solution of (\ref{5}) we have
\begin{equation}
\label{uc}
u \leq C.
\end{equation}
\end{lemma}

\noindent
 \emph{Proof.} This follows from standard blowing-up argument. For the sake of completeness, we include a short proof here. Suppose (\ref{uc}) were not true. Then, there would  exist a sequence of functions $ u_i$ and $L_i \leq L_0$ such that $M_i:=\sup_{x \in \Sigma} u_i (x) \to +\infty$. Without loss of generality, we may assume that $ M_i= u_i (x_i), x_i = (0^{'}, x_{i, N})$.  Indeed, for fixed $i$, there exists a sequence of points $x_i^k$ such that $ u_i (x_i^k) \to \sup_{\bar{\Sigma}} u_i$ as $k \to +\infty$.  Then let $x_i^k=((x_i^k)^{'}, x_{i,N}^k)$ and define $ u_i^k (x)= u_i ( x^{'}+(x_i^k)^{'}, x_N)$. Then, using standard elliptic estimates, one can strike out a sequence $u_i^k$ converging to $u_i^\infty$ as  $k \to +\infty$. We may then replace $u_i$ by $u_i^\infty$ (which we call $u_i$ again). Then, $u_i$ is a solution of (\ref{5}) and there exists $ x_i=(0^{'}, x_N^i)$ such that $ u_i (x_i)=\max_{\bar{\Sigma}} u_i:=M_i$. 

 Then we perform a classical blow up analysis. Set
\begin{equation}
 \ep_i= M_i^{-\frac{p-1}{2}}, v_i (y)= \epsilon_i^{\frac{2}{p-1}} u_i (x_i +\ep_i y).
\end{equation}
Then it is easy to see that $v_i (y)$ satisfies
\[ \Delta v_i - L_i^2 \epsilon_i^2 v_i + v_i^p=0 \ \mbox{in} \ \Sigma_i= \R^{N-1} \times (-\frac{x_{i,N}}{\ep_i}, \frac{L_i-x_{i,N}}{\ep_i}).
\]
Now we extend $v_i$ to $\R^N$ be periodic extension. We still denote the periodic extension as $v_i$. Then, up to extraction of a subsequence,  $ v_i (y) \to v_0(y) $ in $C_{loc}^2 (\R^N)$, and $ v_0$ satisfies the equation
\[ \Delta v_0 + v_0^p=0 \ \mbox{in} \ \R^{N}, v_0(0)=1.\]
However, this is clearly impossible by the  result of Gidas-Spruck
\cite{gs}. \qed

As a corollary, we have

\begin{corollary}
There exists a $L_{*}>0$ such for $L <L_{*}$, $c(L)$ is achieved only by trivial solutions.
\end{corollary}

\noindent \emph{Proof.} Certainly by Sobolev embedding theorem and
Steiner's symmetrization,  a minimizer which is symmetric in $x^{'}$
to $c(L)$ exists. We call it $ u_L$. Now consider the function
$\phi=\frac{\partial u_L}{\partial x_N}$. Then $\phi$ satisfies
\begin{equation}
\label{2.6}
\Delta \phi -L^2 \phi +p u_L^{p-1} \phi=0 \ \mbox{in} \ \Sigma, \ \phi=0 \ \mbox{on} \ \partial \Sigma.
\end{equation}
As $L \to 0$, since $ u_L$ is uniformly bounded by Lemma 2.1 and by the same argument as in  the proof of Lemma 2.1, we conclude that $ \sup_{x \in \Sigma} u_L \to 0$.

On the other hand, since $\phi \in H_0^1 (\Sigma)$, by Poincare's inequality
\begin{equation}
\label{si}
C\int_{\Sigma} |\phi|^2 \leq  \int_{\Sigma} |\nabla \phi|^2
\end{equation}
for some positive constant $C$.

Multiplying (\ref{2.6}) by $\phi$, we obtain using (\ref{si})
\begin{equation}
 p \int_\Sigma u_L^{p-1} \phi^2 =\int_\Sigma |\nabla \phi|^2 +L^2 \int_\Sigma \phi^2 \geq C \int_\Sigma \phi^2.
\end{equation}
 This is impossible if $\phi \not \equiv 0$, for $L$ small enough since $\sup_{x\in \Sigma} u_L \to 0$. Therefore for $L$ small enough, $\frac{\partial u_L }{\partial x_N} \equiv 0$ and $ u_L$ is independent of $x_N$.
\qed

Let us denote by $c^{*} (L)$  the energy level of the trivial solutions, i.e., solutions depending on $x^{'}$ only.  By a simple computation, we see that
\begin{equation}
c^{*} (L)= \gamma_0 L^{2-\frac{(N-1)(p-1)}{p+1}}.
\end{equation}
where $\gamma_0$ is  generic constant (independent of $L$).  Certainly $ c(L) \leq c^{*} (L)$.

\begin{lemma}
As $L\to +\infty$, $ L^{-\frac{2}{p-1}} u_L (L^{-1} y) \to w(y)$ where $w(y)$ is the unique solution of
\begin{equation}
\label{ground}
\left\{\begin{array}{l}
\Delta w -w + w^p=0 \ \mbox{in} \ \R^N,\\
w>0 \ \mbox{in} \ \R^N, w(0)=\max_{y \in \R^N} w(y), w(y) \to 0 \ \mbox{as} \ |y| \to +\infty
\end{array}
\right.
\end{equation}
In fact, $u_L$ has only one local maximum point $P_L$ on $\partial \Sigma$.
\end{lemma}

\noindent \emph{Proof.} By Schwarz spherical rearrangement with
respect to $x^{'}$ and Steiner monotone increasing rearrangement in
$x_N$, after a shift of the origin and a change $x_N$ to $-x_N$ if
needed, we see that  $u_L$ is radially symmetric in $x^{'}$ and
monotone increasing in $ x_N$. Therefore, there exists a unique
point $P_0=(0^{'}, 1)$ where $u_L$ achieves a maximum.  Let $v_L
(y)=L^{-\frac{2}{p-1}} u_L (P_0 +L^{-1} y)$.  Then $ v_L$ satisfies
\[ \Delta v_L -v_L + v_L^p=0 \ \ \mbox{in} \ \Sigma_L, v_L\in H^1 (\Sigma_L) \]
where $\Sigma_L = \R^{N-1} \times (-L, 0)$. Simple computations show that
\begin{equation}
c(L) = (\int_\Sigma u_L^{p+1})^{\frac{p-1}{p+1}} = L^{2-\frac{N(p-1)}{p+1}} (\int_{\Sigma_L} v_L^{p+1})^{\frac{p-1}{p+1}}.
\end{equation}
Since $w$ decays exponentially, we can use a cut-off of $w$ to be a test function and derive that
\begin{equation}
\label{346}
c(L) \leq (\frac{1}{2} \int_{\R^N} w^{p+1} + O(L^{-1}))^{\frac{p-1}{p+1}}
 L^{2-\frac{N(p-1)}{p+1}}
\end{equation}
for $L>>1$. We conclude that  $\int_{\Sigma_L} v_L^{p+1} $ is
bounded and hence the $H^1(\Sigma_L)$ norm  of $v_L$ is bounded.
Consequently,  standard arguments show that as $ L \to +\infty$, $
v_L (y) \to w(y)$. \qed

The following lemma gives part of (1) of Theorem \ref{1.1}.

\begin{lemma}
\label{key}
Let $L_0$ be such that $ c(L_0) < c^{*} (L_0)$. Then for any $L >L_0$, it holds that $c(L)<c^{*} (L)$.
\end{lemma}

\noindent \emph{Proof.} Since $ c(L_0) < c^{*} (L_0)$, $c(L_0)$ is
attained by a nontrivial solution, which we call $ v_0(x^{'}, x_N)$.
Thus, $v_0$ satisfies
\begin{equation}
\Delta v_0 -L_0^2  v_0  + v_0^p=0 \ \mbox{in} \ \Sigma.
\end{equation}
Note that since $v_0$ is nontrivial, we have
\begin{equation}
\int_{\Sigma } v_{0,x_N}^2 >0.
\end{equation}
Now let us consider the following transformation:
\begin{equation}
u(x^{'},x_N)= \lambda^{\frac{2}{p-1}} v_0 (\lambda x^{'}, x_N), \ \mbox{where} \ \lambda=\frac{L}{L_0}.
\end{equation}
A simple computation shows that $ u(x)$ satisfies
\begin{equation}
\Delta u - L^2 u + u^p+ (\lambda^2-1) u_{x_N x_N} =0.
\end{equation}
Hence
\begin{equation}
\int_{\Sigma } (|\nabla u|^2 +L^2 u^2)= \int_{\Sigma } u^{p+1} - (\lambda^2 -1) \int_{\Sigma } |u_{x_N}|^2.
\end{equation}
Now, since $\lambda >1$, we derive the following sequence of inequalities
\begin{eqnarray*}
c(L) & \leq &  \frac{ \int_\Sigma (|\nabla u|^2 + L^2 u^2)}{ (\int_\Sigma u^{p+1})^{\frac{2}{p+1}}}
\\
& < & (\int_\Sigma u^{p+1})^{\frac{p-1}{p+1}} =\lambda^{2 -\frac{(N-1)(p-1)}{p+1}}  (\int_{\Sigma} v_0^{p+1})^{\frac{p-1}{p+1}}
\\
& < &  \lambda^{2 -\frac{(N-1)(p-1)}{p+1}} c^{*} (L_0) =c^{*} (L).
\end{eqnarray*}
\qed

The following inequality may be of independent interest.

\begin{lemma}
Let $u_L$ be a least energy solution to $c(L)$. Then it holds
\begin{equation}
\label{eih}
\int_{\Sigma} [ (|\nabla \varphi |^2 + L^2 \varphi^2 )- p u_L^{p-1} \varphi^2] + (p-1) \frac{ (\int_\Sigma u_L^p \varphi)^2}{ \int_{\Sigma} u_L^{p+1}} \geq 0, \forall \varphi \in H^1 (\Sigma).
\end{equation}
\end{lemma}

\noindent \emph{Proof.} This follows from the variational
characterizations of $u_L$. In fact, let
\[ Q[u]= \frac{ \int_{\Sigma} (|\nabla u|^2+L^2 u^2)}{(\int_{\Sigma} u^{p+1})^{\frac{2}{p+1}}}
\]
and $ \rho (t)= Q[u_L+ t \varphi]$ for any $\varphi \in H^1 (\Sigma)$. Since $\rho (t) \geq \rho (0)$ for any $t$, $\rho^{'} (0)=0,\rho^{''} (0) \geq 0$.  Note that
\begin{eqnarray*}
  & \ \ & \int_\Sigma (|\nabla (u_L+t \varphi)|^2 +L^2(u_L+t \varphi)^2 )\\
  &=& \int_\Sigma (|\nabla u_L|^2 +L^2 u_L^2)
  +2 t \int_\Sigma (\nabla u_L \nabla \varphi+L^2 u_L \varphi) + t^2 \int_\Sigma (|\nabla \varphi|^2+L^2 \varphi^2)
  \\
  &=& \int_\Sigma u_L^{p+1}
  +2 t \int_\Sigma u_L^p \varphi + t^2 \int_\Sigma (|\nabla \varphi|^2+L^2 \varphi^2)
  \\
 & \ \ & \int_{\Sigma} (u_L+t \varphi)^{p+1}= \int_{\Sigma} u_L^{p+1} + (p+1) t \int_\Sigma u_L^p \varphi + \frac{(p+1)p}{2}\int_\Sigma u_L^{p-1}\varphi^2 + {\mathcal O} (t^2)
\end{eqnarray*}
Then that $\rho^{''} (0)\geq 0$ is equivalent to
\begin{equation}
\label{eih1}
\int_{\Sigma} [ (|\nabla \varphi |^2 + L^2 \varphi^2 )- p u_L^{p-1} \varphi^2] + (p-1) \frac{ (\int_\Sigma u_L^p \varphi)^2}{ \int_{\Sigma} u_L^{p+1}} \geq 0.
\end{equation}
\qed

The next result is a corollary of  inequality (\ref{eih}).
\begin{lemma}
\label{key1}
Let $u_L$ be a least energy solution of $c(L)$ and $\lambda_2 (u_L)$ be the second eigenvalue of
\begin{equation}
\Delta \phi-L^2 \phi+pu_L^{p-1} \phi+\lambda \phi=0 \ \mbox{in} \ \Sigma, \ \frac{\partial \phi}{\partial \nu}=0 \ \mbox{on} \ \partial \Sigma.
\end{equation}
Then, necessarily
\begin{equation}
\lambda_2 (u_L) \geq 0.
\end{equation}
\end{lemma}

\noindent \emph{Proof.} Recall that by the Courant-Fisher-Weyl
formula, one has
\begin{equation}
\label{345}
 \lambda_2 = \max_{ dim (V)=1, V \in H^1 (\Sigma)} \inf_{\varphi \perp V}\int_{\Sigma} [ (|\nabla \varphi |^2 + L^2 \varphi^2 )- p u_L^{p-1} \varphi^2].
 \end{equation}
 Then by choosing $V=\mbox{span} \{ u_L^p\}$ in (\ref{345}) and using (\ref{eih}), we derive that $\lambda_2\geq 0$.
\qed

\vskip 0.5cm

\noindent
{\bf Completion of proof of (1) of Theorem \ref{1.1}:} Let
\[ L_{*}= \sup \{ L |  c(l) =c^{*} (l) \ \mbox{for}  \ l \in (0, L) \}.\]
 By Corollary 2.2, we see that $0<L_{*}$. Now from  Lemma 2.3, it follows that $L_{*}<+\infty$. Indeed, for $L$ large, we have by (\ref{346}) and Lemma 2.3, $c(L) \sim L^{2-\frac{N(p-1)}{p+1}}$, while $c^{*}(L) \sim L^{2-\frac{(N-1)(p-1)}{p+1}}$.  Certainly, for $ l \leq L_{*}$, $c(l)=c^{*} (l)$ and $u_l$ is trivial. By Lemma \ref{key}, $ c(L) <c^{*} (L)$ for $L >L_{*}$.

We now claim that $L^{*}=\frac{\pi}{\sqrt{\lambda_1}}$.  In fact, by separation of variables, the second eigenvalue of the following eigenvalue problem
\begin{equation}
\Delta \phi-L^2 \phi+pw_0^{p-1} \phi+\lambda \phi=0 \ \mbox{in} \ \Sigma, \ \frac{\partial \phi}{\partial \nu}=0 \ \mbox{on} \ \partial \Sigma
\end{equation}
is
\begin{equation}
\label{367}
\lambda_2  =\pi^2 -L^2\lambda_1  \geq 0.
\end{equation}

By Lemma \ref{key1}, $L^{*} \leq \frac{\pi}{\sqrt{\lambda_1}}$. On the other hand, if $L^{*}>\frac{\pi}{\sqrt{\lambda_1}}$, then by Lemma \ref{key}, $w_{0}$ is a minimizer of $c_L$ when $L$ is close to $\frac{\pi}{\sqrt{\lambda_1}}$. But it is easy to see that $w_{0}$ loses its stability exactly at $L=\frac{\pi}{\sqrt{\lambda_1}}$ by (\ref{367}).

\section{ The subcritical Case: the proof of (2) of Theorem \ref{1.1} }
\setcounter{equation}{0}

In this section, we prove the uniqueness and nondegeneracy of the least energy solutions for $L \to +\infty$.

First we recall

\begin{lemma}
\label{l3.2}
Let $w$ be the least energy solution of
\begin{equation}
\label{ground1}
\Delta w- w + w^p=0, w>0,  w \in H^1 (\R^{N}).
\end{equation}
Then $w$ is nondegenerate, i.e.,
\begin{equation}
\mbox{Ker} \ (\Delta -1 + pw^{p-1})= \mbox{span} \ \{ \frac{\partial w}{\partial y_1}, ..., \frac{\partial w}{\partial y_N} \}
\end{equation}
\end{lemma}

\noindent \emph{Proof.} The result is well-known. By the clasical
result of Gidas-Ni-Nirenberg \cite{gnn}, $w$ is radially  symmetric.
The nondegeneracy follows from the uniqueness result of Kwong
\cite{k}. See Lemma A.3 of \cite{nt2}. Here we include a short and
self-contained  new proof of nondegeneracy using only property of
least energy. This part is of independent interest, but it is
restricted to the power nonlinearity.

 Let $w=w(r)$ be a radial least energy solution of (\ref{ground}). By the same proof as in Lemma \ref{key1}, $ \lambda_{2,r} (w) \geq 0$, where $\lambda_{2,r}$ denotes the second eigenvalue in the radial class. It remains to show that $\lambda_{2,r} (w) >0$. Suppose $\lambda_{2,r} =0$ and let $\phi (r)$ be the corresponding eigenfunction, i.e.
\begin{equation}
\label{phi1}
\Delta \phi -\phi + pw^{p-1} \phi=0, \phi=\phi (r)\in H^1 (\R^N).
\end{equation}
Then the characterization of the second eigenfunction implies that $\phi$ changes sign once. So we may assume that $ \phi <0$ for $ r< r_0$ and $ \phi >0$ for $r>r_0$. Now as in Kwong-Zhang \cite{kz}
 we consider the function
\begin{equation}
\eta (r)= rw^{'} -\beta w.
\end{equation}
Then $\eta$ satisfies
\begin{equation}
\label{phi2}
\Delta \eta-\eta+ pw^{p-1} \eta=  2w - (2+\beta (p-1)) w^p.
\end{equation}
We choose $\beta$ such that $ 1= (1+\frac{\beta (p-1)}{2}) w^{p-1} (r_0)$, hence  $2w - (2+\beta (p-1)) w^p <0$ for $r <r_0$ and $  2w - (2+\beta (p-1)) w^p >0$ for $r>r_0$.

Multiplying (\ref{phi1}) by $ \eta$ and (\ref{phi2}) by $\eta$, we arrive at
\begin{equation}
\int_{\R^N} \phi  (2w - (2+\beta (p-1)) w^p ) =0
\end{equation}
which is impossible by the property of $\phi$. Thus $\phi \equiv 0$
and this completes the proof. \qed

Let us now prove the nondegeneracy of the least energy solution when $L$ is large. By the rescaling $ w_L= L^{-\frac{2}{p-1}} u_L ((0,1)+L^{-1} y)$,
 it is enough to show that the only solutions to
\begin{equation}
\label{non1}
\Delta \phi -\phi + pw_{L}^{p-1} \phi=0, \phi \in H^1 (\Sigma_L)
\end{equation}
are $\frac{\partial w_L}{\partial y_j}, j=1,..., N-1$. Here $\Sigma_L= \R^{N-1} \times (-L, 0)$. We may assume that
\begin{equation}
\label{non2}
\int_{\Sigma_L} \phi \frac{\partial w_L}{\partial y_j} =0, j=1,..., N-1.
\end{equation}

Suppose there is a nonzero solution $\phi$ to (\ref{non1})-(\ref{non2}). We may assume that $ \| \phi \|_{L^\infty} =1$.  Since $w_L \to w(y)$ as $L\to \infty$, we conclude that $p w_L^{p-1} <\frac{1}{2}$ for $y \in B_R (0) \cap \Sigma_L$. Thus $ |\phi (y)| \leq \max_{y \in B_R (0)} |\phi| e^{-\frac{1}{\sqrt{2}} (|y|-R)}$. So the maximum point of $ |\phi|$ must occur in $B_{2R} (0)$.  Letting $L \to +\infty$, we have that $\phi \to \phi_{\infty}$ which satisfies
\begin{equation}
\label{non3}
\Delta \phi_{\infty} -\phi_{\infty} + pw^{p-1} \phi_\infty=0, \int_{\R^N} \phi_{\infty} \frac{\partial w}{\partial y_j}=0, j=1, ..., N-1.
\end{equation}
By Lemma \ref{l3.2}, $ \phi_\infty= c \frac{\partial w}{\partial y_N}= c \frac{w^{'}}{r} y_N$. We have seen that $\phi_\infty$ attains its maximum at some finite point. This is impossible since $ \frac{\partial^2 w }{\partial^2 y_N} \not = 0$.

The proof of uniqueness of least energy solution when $L$ is large is similar to that of nondegeneracy. In fact, suppose that there are two least energy solutions  $u_{L}$ and $ u_{L}^{'}$ to (\ref{5}). We may assume that both $ u_{L}$ and $u_{L}^{'}$ attain their maximum at $(0,1)$. Suppose that $ u_{L} \not \equiv u_{L}^{'}$. Then letting $ w_L (y):= L^{-\frac{2}{p-1}} u_L ((0,1)+L^{-1} y), w_L^{'} (y):= L^{-\frac{2}{p-1}} u_L^{'} ((0,1)+L^{-1} y) ), \phi_L (y):= w_L (y)- w_L^{'} (y)$, we see that $\phi_L$ satisfies
\begin{equation}
\label{non1-1}
\Delta \phi -\phi +  V(y) \phi=0, \phi \in H^1 (\Sigma_L)
\end{equation}
 where $V(y)= \frac{ w_L^p- (w_L^{'})^p}{w_L -w_L^{'}}$. Since $V(y) \to pw^{p-1} (y)$ as $L \to +\infty$ and $\nabla \phi_L (0)=0$, the rest of the proof is exactly the same as before. We omit the details.

\section{ The critical exponent Case: the proof of (3) of Theorem \ref{1.1} }
\setcounter{equation}{0}

In this section, we assume that $p=\frac{N+2}{N-2}$. We consider two cases: $L$ is small and $L$ is large

\subsection{Critical exponent case I: $L$ small}

 It is well-known (see \cite{gs}) that the solutions to the following problem
\begin{equation}
\label{critical}
\Delta u + u^{\frac{N+2}{N-2}}=0 \ \mbox{in} \ \R^N, u>0
\end{equation}
are given by
\begin{equation}
U_{\epsilon, a} = c_N (\frac{\ep}{\ep^2+|x-a|^2})^{\frac{N-2}{2}}
\end{equation}
for some $\epsilon >0 $ and $ a \in \R^N$.

Let
\begin{equation}
\label{S}
S= \frac{\int_{\R^N} |\nabla U_{1,0}|^2}{(\int_{\R^N} U_{1,0}^{\frac{2N}{N-2}})^{\frac{N-2}{N}}}, \ S_{\frac{1}{2}}= (\frac{1}{2})^{\frac{2}{N}} S.
\end{equation}

We have the following lemma, whose proof follows from classical
``concentration-\break compactness'' principle  of P.L. Lions
\cite{pl1}, \cite{pl2}.

\begin{lemma}
\label{cc}
Let $p=\frac{N+2}{N-2}$. If
\begin{equation}
\label{c1}
c(L) <S_{\frac{1}{2}}
\end{equation}
then $c(L)$ is attained.
\end{lemma}

As a Corollary, we have

\begin{corollary}
\label{c4.2}
Let $N \geq 4$. Then for $L$ sufficiently small, $c(L)$ is attained.
\end{corollary}

\noindent \emph{Proof.} We just need to verify (\ref{c1}) for $L$
small. Now we compute
\begin{eqnarray}
 & & \int_\Sigma U_{\ep, 0}^{\frac{2N}{N-2}}\\
 &=&  c_N^{\frac{2N}{N-2}} \int_\Sigma (\frac{\ep}{\ep+|x|^2})^N dx\nonumber \\
&=& c_N^{\frac{2N}{N-2}} \int_0^{\frac{1}{\ep}} (\int_{\R^{N-1}} \frac{1}{ (1+t^2 +|y^{'}|^2)^N} d y^{'} dt\nonumber \\
& = & c_N^{\frac{2N}{N-2}} \int_0^{\frac{1}{\ep}} (1+t^2)^{-\frac{N+1}{2}} dt (\int_{\R^{N-1}} \frac{1}{ (1 +|y^{'}|^2)^N} d y^{'} \nonumber \\
& = & c_N^{\frac{2N}{N-2}}  \Bigg[\int_{\R^{N}} \frac{1}{ (1+|y|^2)^N} d y - \frac{\ep^N }{N} \int_{\R^{N-1}} \frac{1}{ (1+|y^{'}|^2)^N} d y^{'} + o(\ep^N) \Bigg]
\label{cal1}.
\end{eqnarray}

Similarly
\begin{equation}
\label{cal2}
 \int_\Sigma |\nabla U_{\ep, 0}|^2 = c_N^2 (N-2)  \int_{R^N_{+}}  |\nabla U_{1,0}|^2 dy
\end{equation}
\[-c_N^2 \ep^{N-2} \int_{\R^{N-1}} \frac{1}{(1+|y^{'}|^2)^{N-1}} d y^{'} +  c_N^2 \frac{\ep^N  (N-2)}{N} \int_{\R^{N-1}} \frac{1}{ (1+|y^{'}|^2)^N} d y^{'} + O(\ep^N). \]

On the other hand,  for $N\geq 5$,
\begin{equation}
\label{cal3}
 \int_\Sigma U_{\ep, 0}^2=
c_N^2 \int_\Sigma (\frac{\ep}{\ep^2+|x|^2})^{N-2}= c_N^2 \ep^2 \int_{\R_{+}^N} \frac{1}{ (1+|y|^2)^{N-2}} dy
\end{equation}
Thus for $N\geq 5$
\[
c(L) \leq \frac{ \int_\Sigma (|\nabla U_{\ep, 0}|^2 +L^2 U_{\ep, 0}^2)}{ (\int_{\Sigma} U_{\ep, 0}^{p+1})^{\frac{N-2}{N}} } = \frac{ A_0 -B_0 \ep^{N-2} +L^2 C_0 \ep^2 +O(\ep^N) }{ (D_0 +O( \ep^N) )^{\frac{N-2}{N}} }
<\frac{A_0}{D_0^{\frac{2}{p+1}}}= S_{\frac{1}{2}}
\]
if $L$ is small and $\ep$ is small. Here
\begin{equation}
A_0= \int_{\R^N_{+}} |\nabla U_{1,0}|^2, B_0= c_N^2 \int_{\R^{N-1}} \frac{1}{ (1+|y^{'}|^2)^{N-1}} d y^{'},  C_0=  \int_{\R_{+}^N} U_{1,0}^2, D_0= \int_{\R_{+}^N} U_{1,0}^{\frac{2N}{N-2}}
\end{equation}

Applying Lemma (\ref{cc}),  for $L$ small, $ c(L)$ is attained (possibly by a trivial solution).

 For $N=4$, we have
\begin{equation}
\label{cal4}
 \int_\Sigma U_{\ep, 0}^2= c_N^2 \ep^2 \log \frac{1}{\ep} \int_{\R^{N-1}} \frac{1}{ (1+|y^{'}|^2)^2} d y^{'}
\end{equation}
Similar arguments as before show that $c(L) < S_{\frac{1}{2}}$ for $L$ small.
\qed

\bigskip
\noindent
{\bf Remark:} In fact, Corollary  4.2 is also true for $N\geq 3$. Another  proof is to use the inequality $ c(L) \leq c^{*}(L)$.

\subsection{Critical exponent case II: $L$ large}

The main theorem in this section is the following

\begin{theorem}
\label{4.1}
For $L$ large, $ c(L)$ is not attained.
\end{theorem}

We first assume that $N\geq 5$. Later on
 we will show how one can modify the arguments to deal with the case of $N=4$.

 We prove it by contradiction.  Suppose that  $ c(L)$ is attained by some $u_L$ for a sequence of $L=L_i \to +\infty$.  Note that the computations of Corollary \ref{c4.2} show that
\begin{equation}
\label{clbd}
c(L) \leq S_{\frac{1}{2}}.
\end{equation}
We claim that $ c(L)$ is attained for all $L $. In fact, let $ L  < L_{i}$ for some $L_i$. Then we have
\begin{equation}
c(L) \leq  \frac{ \int_\Sigma (|\nabla u_{L_i} |^2 + L^2 u_{L_i}^2)}{ (\int_\Sigma u_{L_i}^{p+1})^{\frac{2}{p+1}}} <c(L_i) \leq S_{\frac{1}{2}}.
\end{equation}
By Lemma \ref{cc}, $c(L)$ is attained.

By Steiner symmetrization, we may assume that $ u_L$ is symmetric in $ x^{'}$ and increasing in $x_N$.

 We rescale $ \tilde{u}_L= L^{-\frac{2}{p-1}} u_L (L^{-1} y)$ and obtain that
\[ c(L) = \frac{ \int_{\Sigma_L} ( |\nabla \tilde{u}_L|^2 + \tilde{u}_L^2)}{(\int_{\Sigma_L} \tilde{u}_L^{p+1})^{\frac{2}{p+1}}} \leq S_{\frac{1}{2}}
\]
where $\Sigma_L= \R^{N-1} \times (0, L)$. This implies that
\begin{equation}
\int_{\Sigma_L} (|\nabla \tilde{u}_L|^2 +\tilde{u}_L^2) \leq C, \int_{\Sigma_L} \tilde{u}_L^{p+1} \leq C.
\end{equation}

We claim that $\tilde{u}_L$ must blow up. If not, by taking a subsequence of $L$ and extending $\tilde{u}_L$ to $ \R^{N-1} \times (-L, L)$, we see that $v_L$ converges to a positive solution of
\[ \Delta v-v + v^{\frac{N+2}{N-2}}=0, v \in H^1 (\R^N)\]
 which is impossible. In fact we have that
\[ \frac{1}{ \tilde{u}_L (0)} \tilde{u}_L ( (\tilde{u}_L (0))^{-\frac{p-1}{2}} (y)
 \to U_{1,0} (y) \]
in $C_{loc}^2 (\R_{+}^N)$ as $L \to \infty$.

Let $\ep$ be such that
\begin{equation}
\label{epdef}
  \tilde{u}_L (0)= \ep^{-\frac{N-2}{2}}
\end{equation}
and set
\begin{equation}
v_\ep (y)= \ep^{-\frac{N-2}{2}} \tilde{u}_L (\ep y).
\end{equation}

Then it is easy to see that $v_\ep$ satisfies
\begin{equation}
\label{vepeqn}
\Delta v_\ep -\ep^2 v_\ep + v_\ep^p=0 \ \mbox{in} \ \Sigma_{\frac{L}{\ep}}, \frac{\partial v_\ep}{\partial \nu}=0 \ \mbox{on} \ \Sigma_{\frac{L}{\ep}}
\end{equation}
and that $v_\ep  (y) \to U_{1,0} (y) $ in $C_{loc}^2 (\R_{+}^N)$.

We now require the following crucial estimate.

\begin{lemma}
\begin{equation}
v_\ep (y) \leq \frac{C}{ ( 1+|y|^2)^{\frac{N-2}{2}}}.
\end{equation}
\end{lemma}

\noindent \emph{Proof.} The derivation of this estimate   follows
exactly the proof of Theorem 2.1 of \cite{gl}. In fact, our
situation is simpler as there is no need to straighten  the boundary
(see \cite{gl}). \qed

Let  $V_{\ep}$    be the unique solution of the following linear problem
\begin{equation}
\Delta V_{\ep} - \ep^2 V_{\ep } + U_{1, 0}^{\frac{N+2}{N-2}}=0 \ \mbox{in} \ \Sigma_{\frac{L}{\ep}}, \ \
\frac{\partial V_{\ep}}{\partial \nu} =0 \ \mbox{on} \ \partial \Sigma_{\frac{L}{\ep}}.
\end{equation}
Now we decompose
\[ V_{\ep}= U_{1, 0}-\varphi_\ep\]
Then we have
\begin{equation}
\label{projection1}
\left\{\begin{array}{l}
\Delta \varphi_{\ep} - \ep^2  \varphi_{\ep} +\ep^2  U_{1, 0}=0 \ \mbox{in} \ \Sigma_{\frac{L}{\ep}}, \\
\frac{\partial \varphi_\ep}{\partial \nu} =\frac{\partial U_{1, 0}}{\partial \nu}  \ \mbox{on} \ \partial \Sigma_{\frac{L}{\ep}}.
\end{array}
\right.
\end{equation}
Let us set
\begin{equation}
\varphi_\ep=  \ep^2 \varphi_0 (y) +\varphi_{\ep, 1}
\end{equation}
where $\varphi_0$ is the unique solution of the following problem
\begin{equation}
\label{var}
\Delta \varphi_0 + U_{1,0} (y)=0 \ \mbox{in} \ \R^N, \varphi_0 (y)=\varphi_0 (|y|), \varphi_0 \to 0 \ \mbox{as} \ |y| \to +\infty
\end{equation}
Note that since
\[ U_{1,0} (y) \leq \frac{C}{ (1+|y|)^{N-2}}
\]
where $N-2 >2$, there exists a unique solution to (\ref{var}).

Then $\varphi_{\ep,1}$ satisfies
\begin{equation}
\label{projection2}
\left\{\begin{array}{l}
\Delta \varphi_{\ep,1} -   \ep^2 \varphi_{\ep,1} -  \ep^4 \varphi_0 (y) =0 \ \mbox{in} \ \Sigma_{\frac{L}{\ep}}, \\
\frac{\partial \varphi_{\ep, 1}}{\partial \nu} =\frac{\partial }{\partial \nu} [U_{1,0 } -(\ep^2) \varphi_0 ]  \ \mbox{on} \ \partial \Sigma_{\frac{L}{\ep}}.
\end{array}
\right.
\end{equation}

We claim that

\begin{lemma}
\label{4.3}
\begin{equation}
|\varphi_{\ep,1} | \leq  C \ep^{N-2} +o( \ep^2) \frac{1 }{ (1+|y|^2)^{\frac{N-2}{2}}}.
\end{equation}
\end{lemma}

Now we let
\begin{equation}
  v_\ep ( y)=  V_{\ep}  (y) +\phi_\ep (y)
\end{equation}
Then $\phi_\ep$ satisfies
\begin{equation}
\label{3.25}
\left\{\begin{array}{l}
 \Delta \phi_\ep - \ep^2 \phi_\ep +  ( V_{\ep} +\phi_\ep )^{\frac{N+2}{N-2}} - U_{1, 0}^{\frac{N+2}{N-2}}=0 \ \mbox{in} \ \Sigma_{\frac{L}{\ep}},\\
\frac{\partial \phi_\ep}{\partial \nu}=0 \ \mbox{on} \ \ \partial \Sigma_{\frac{L}{\ep}}
\end{array}
\right.
\end{equation}
We also claim that
\begin{lemma}
\label{4.4}
  \begin{equation}
| \phi_\ep (y) | \leq C  \ep^2 \frac{1}{ (1+|y|^2)^{\frac{N-2}{2}}}.
\end{equation}
\end{lemma}

Postponing the proofs of Lemma (\ref{4.3}) and Lemma (\ref{4.4}) to the appendix, we can conclude that  we have a contradiction to (\ref{clbd}) by establishing the following:
\begin{equation}
\label{clld}
 c(L) >S_{\frac{1}{2}}.
\end{equation}

First we note that
\begin{equation}
\label{clnn}
 c(L)= \frac{ \int_{\Sigma_{\frac{L}{\ep}}} (|\nabla v_\ep|^2 + \ep^2 v_\ep^2)}{(\int_{\Sigma_{\frac{L}{\ep}}} v_\ep^{p+1})^{\frac{2}{p+1}}}.
\end{equation}

Then we have
\begin{eqnarray*}
 \int_{\Sigma_{\frac{L}{\ep}}} (|\nabla v_\ep|^2 + \ep^2 v_\ep^2) &=&  \int_{\Sigma_{\frac{L}{\ep}}} (|\nabla (V_\ep +\phi_\ep)|^2 + \ep^2 (V_\ep +\phi_\ep)^2)
\nonumber \\
& = & \int_{\Sigma_{\frac{L}{\ep}}} (|\nabla V_\ep |^2 + \ep^2 V_\ep^2)
+2 \int_{\Sigma_{\frac{L}{\ep}}} (\nabla V_\ep \nabla \phi_\ep + \ep^2 V_\ep \phi_\ep)
\nonumber \\
&  + & \int_{\Sigma_{\frac{L}{\ep}}} (|\nabla \phi_\ep |^2 + \ep^2 \phi_\ep^2) \nonumber \\
& = & I_1+2 I_2+I_3
\end{eqnarray*}
where $I_1, I_2$ and $I_3$ are defined by the three terms at the last equality.

Quantity $I_1$  can be computed as follows:
\[ I_1=\int_{\Sigma_{\frac{L}{\ep}}} U_{1,0}^p V_{\ep}
 = \int_{\Sigma_{\frac{L}{\ep}}} U_{1,0}^p (U_{1,0} -\varphi_{\ep})  \]
\begin{equation}
\label{i1}
= \int_{\R_{+}^N } U_{1,0}^{p+1} - \ep^2 \int_{\R_{+}^N} U_{1,0}^2  + o(\ep^2)
\end{equation}
where we have used
\begin{equation}
  \int_{\R_{+}^N} U_{1,0}^p \varphi_0 = \int_{\R_{+}^N} U_{1,0}^2 >0.
\end{equation}

 For  the quantity $I_2$, from the equation for $ V_\ep$, it follows that
\begin{equation}
\label{i2}
 I_2= \int_{\Sigma_{L}{\ep}} U_{1,0}^p \phi_\ep = O(\ep^2).
\end{equation}
By Lemma \ref{4.3}, we have
\[ I_3=  \int_{\Sigma_{\frac{L}{\ep}}} (|\nabla (\phi_\ep)|^2 + \ep^2 (\phi_\ep)^2) = o(\ep^2).\]

So
\begin{equation}
\label{cl11}
 \int_{\Sigma_{\frac{L}{\ep}}} (|\nabla v_\ep|^2 + \ep^2 v_\ep^2)  \geq  \int_{\R_{+}^N} |\nabla U_{1,0}|^2 - \ep^2  \int_{\R_{+}^N} U_{1,0}^2  + 2 I_2 + o(\ep^2).
\end{equation}

Next, using Lemma \ref{4.3} again, we obtain that
\begin{eqnarray}
 \int_{\Sigma_{\frac{L}{\ep}}} v_\ep^{p+1}&=&  \int_{\Sigma_{\frac{L}{\ep}}} (V_\ep+\phi_\ep)^{p+1}
\nonumber \\
& = &  \int_{\Sigma_{\frac{L}{\ep}}} V_\ep^{p+1} + (p+1)   \int_{\Sigma_{\frac{L}{\ep}}} V_\ep^p \phi_\ep + o(\ep^2)
\nonumber \\
 &= &  \int_{\Sigma_{\frac{L}{\ep}}} V_\ep^{p+1} + (p+1)   \int_{\Sigma_{\frac{L}{\ep}}} U_{1,0}^p \phi_\ep + o(\ep^2) \nonumber \\
& = & \int_{\R_{+}^N} U_{1,0}^{p+1} - (p+1) \ep^2 \int_{\R_{+}^N} U_{1,0}^2  + I_2 + o(\ep^2).
\label{cl2}
\end{eqnarray}

Combining (\ref{cl11}) and (\ref{cl2}), we obtain that
\begin{eqnarray}
\label{cl3}
c(L) &\geq & \frac{  \int_{\R_{+}^N} |\nabla U_{1,0}|^2 -  \ep^2 \int_{\R_{+}^N} U_{1,0}^2 + I_2}{ ( \int_{\R_{+}^N} U_{1,0}^{p+1} - (p+1) \ep^2 \int_{\R_{+}^N} U_{1,0}^2 + (p+1) I_2 + o(\ep^2))^{\frac{2}{p+1}}}\\
& \geq & S_{\frac{1}{2}} + \ep^2 \int_{\R_{+}^N} U_{1,0}^2 +o(\ep^2) + O(|I_2|^2) >S_{\frac{1}{2}} \nonumber
\end{eqnarray}
which proves (\ref{clld}). \qed

Finally, when $N=4$, we have to replace $ \varphi_0 (y)$ be the following function
\[ \varphi_0 (y)= \log \frac{1}{1+|y|^2} + (\Delta)^{-1} (\frac{1}{ (1+|y|^2)^2})\]
and $\ep^2$ by $ \ep^{2} \log \frac{1}{\ep}$. The rest of the proof is unchanged.

\subsection{Completion of proof of (3) of Theorem 1.1}

Let $p=\frac{N+2}{N-2}$ and $N\geq 4$.  By Corollary \ref{c4.2}, $c(L)$ is attained for $L$ small. Set
\begin{equation}
L_{**}=\sup \{ L|  c(L) \ \mbox{is attained for} \ l \in (0, L) \}.
\end{equation}

By Corollary \ref{c4.2} and Theorem \ref{4.1}, $0<L_{**} <+\infty$.

We claim that $c(L)$ is also attained at $L=L_{**}$. In fact, if not, let $u_{L}, L<L_{**}$ be the minimizers of $c(L)$. Then as $L \to L_{**}$, $ u_{L}$ must blow up. But similar arguments as in Theorem \ref{4.1} shows that this is impossible.

Now we show that for $L>L_{**}$, $c(L)$ is not attained. In fact, suppose $ c(L)$ is attained for some $L_0 >L_{**}$. then certainly, $c(L_0) \leq S_{\frac{1}{2}}$. Let the minimizer of $ c(L_0)$ be $u_{L_0}$. Then
\begin{equation}
c(L) <c(L_0) \leq S_{\frac{1}{2}}, \ \mbox{for} \ L <L_0.
\end{equation}
By Lemma \ref{cc}, $c(L)$ is attained if $L < L_0$, which is a contradiction with  the definition of $ L_{**}$.

Finally for $L$ small, we have $ c(L) <c(\frac{L_{**}}{2}) \leq S_{\frac{1}{2}}$. By the same analysis as in Lemma \ref{l00} the minimizer is uniformly bounded. Hence for $L$ small $ c(L)$ is archived by  trivial solution.  Now similar proofs as in Lemma \ref{key} yield another constant $L_{*} <L_{**}$ such that for $L \leq L_{*}$, $c(L)$ is achieved by trivial solution and  for $L \in (L_{*}, L_{**}]$, $c(L)$ is attained by a nontrivial constant.

Thus, (3) of Theorem \ref{1.1} is  proved. \qed

\bigskip

\noindent
{\bf  Appendix A: Proof of Lemma \ref{4.3} and Lemma \ref{4.4}.}

\bigskip

To prove Lemma \ref{4.3} and \ref{4.4}, we introduce two weighted $L^\infty$ spaces. For $f$ a function in $\Sigma_{\frac{L}{\ep}}$, we define the following weighted $L^\infty$-norms
\[ \| f\|_{*}= \sup_{x \in \Sigma_{\frac{L}{\ep}}} \Big\vert(1+|y|^2)^{\frac{N-2}{2}}f(x)\Big\vert\]
and
\[
\Vert f\Vert_{**}=\sup_{x\in\Sigma_{\frac{L}{\ep}}}\Big\vert(1+|y|^2)^{\frac{N-1}{2}} f(x)\Big\vert.
\]

\noindent
{\it {\bf Lemma A.}
Let $ f \in L^\infty (\Sigma_{\frac{L}{\ep}})$ be such that
\[ \| f\|_{**} <+\infty\]
 and $ u$ satisfy
\[ -\Delta u + \ep^2 u = f \ \ \mbox{in} \ \Sigma_{\frac{L}{\ep}}, \ \frac{\partial u}{\partial \nu} =0 \ \mbox{on} \ \partial \Sigma_{\frac{L}{\ep}}.\]
Then we have
\begin{equation}
\label{uupper}
| u(x)| \leq  \int_{\Sigma_{\frac{L}{\ep}}} \frac{C}{|x-y|^{N-2} } |f(y) | dy,
\end{equation}
where $C$ is independent of $L \geq L_0$.
As a consequence, we have
\begin{equation}
\| u \|_{*} \leq C \| f \|_{**}.
\end{equation}
}

To prove Lemma \ref{4.3}, we decompose $\varphi_{\ep,1}$ into two parts:
\[ \varphi_{\ep, 1}=\varphi_{\ep, 1}^1+\varphi_{\ep, 1}^2\]
where $\varphi_{\ep, 1}$ satisfies
\[ \Delta \varphi_{\ep, 1}^1- \ep^2 \varphi_{\ep, 1}^1=0 \ \mbox{in} \ \Sigma_{\frac{L}{\ep}}, \frac{\partial \varphi_{\ep,1}^1}{\partial \nu}=\frac{\partial}{\partial \nu} (U_{1, 0} - \ep^2  \varphi_0 (y)) \ \mbox{on} \ \partial \Sigma_{\frac{L}{\ep}},
\]
and  $\varphi_{\ep, 1}^2$ satisfies
\[ \Delta \varphi_{\ep, 1}^2- \ep^2 \varphi_{\ep, 1}^2-  \ep^2  \varphi_0 (y)=0 \ \mbox{in} \ \Sigma_{L}, \frac{\partial \varphi_{\ep,1}^2}{\partial \nu}=0 \ \mbox{on} \ \partial \Sigma_{\frac{L}{\ep}}
\]

The first part can be estimated by asymptotic analysis while the second part follows from comparison principle.

To prove Lemma \ref{4.4}, we note that $\phi_\ep$ satisfies
\[ \Delta \phi_\ep -\ep^2 \phi_\ep + p V_{\ep}^{p-1} \phi_\ep  +E_\ep + N_\ep [\phi_\ep]=0\]
where
\[E_\ep=  V_{\ep}^{p}- U_{1, 0}^p, \]
\[ N_\ep [\phi_\ep]= ( V_{\ep} +\phi_\ep)^p-  V_{\ep}^{p} - p   V_{\ep}^{p-1} \phi_\ep. \]

We argue by contradiction. Let $\tilde{\phi}_\ep=\frac{\phi_\ep}{ \ep^2}$. We assume that
\begin{equation}
\label{A1}
\|\tilde{\phi}_\ep\|_{*} \to +\infty.
\end{equation}
Set
\begin{equation}
\label{A2}
 \Phi_\ep = \frac{\tilde{\phi}_\ep}{\| \tilde{\phi}_\ep\|_{*} }.
\end{equation}

Then it is easy to see that $\Phi_\ep$ satisfies
\[ \Delta \Phi_\ep - \ep^2 \Phi_\ep +  p  V_{\ep}^{p-1} \Phi_\ep
 +  ( \| \tilde{\phi}_\ep\|_{*} )^{-1}
\ep^{-2} E_\ep+  ( \| \tilde{\phi}_\ep\|_{*} )^{-1} \ep^{-2} N_\ep=0.\]

As $\ep \to 0$, $\Phi_\ep \to \Phi_0$ where $\Phi_0$ satisfies
\[ \Delta \phi_0 + p U_{1,0}^{p-1} \phi_0=0 \ \mbox{in} \ \R^N.\]
It is well-known (see \cite{rey}) that $ \Phi_0= a_0 \frac{\partial U_{\lambda, 0}}{\partial \lambda}|_{\lambda=1}+\sum_{j=1}^{N-1} a_j \frac{\partial U_{1,0}}{\partial y_j}$ for some constants $a_j, j=0,1, ..., N-1$.

Now since both $v_\ep$ and $V_{\ep}$  are symmetric in $x^{'}$, we see that $\frac{\partial}{\partial y_j} \Phi_\ep (0)=0, j=1,..., N-1$. We also have that $ \Phi_\ep (0)= \frac{ v_\ep (0) - V_{\ep} (0)}{  \ep^2 \|\tilde{\phi}_\ep\|_{*}}= o(1)$ and hence $\Phi_0(0)=0$. This together with $\frac{\partial}{\partial y_j} \phi_0 (0)=0$ will force $a_j=0, j=0,1, ..., N-1$ and hence $\Phi_0 =0$.

On the other hand, from the equation for $\Phi_\ep$, we have that
\[ \|   V_{\ep}^{p-1} \Phi_\ep \|_{**} = o(1), \| \ep^{-2} E_\ep\|_{**}=O(1),\]
and by Lemma A, we then arrive at
\[ \|\Phi_\ep\|_{*}= o(1)\]
A contradiction to (\ref{A2})!

So Lemma \ref{4.3} is proved. \qed

It remains to prove Lemma A. By a scaling, we may assume that $\ep=1$. Then
\[ u(x)= \int_{\Sigma_L} G_L (x, y) f(y) dy\]
where $G_L (x, y)$ is the Green's function
\[ \Delta G_L -G_L +\delta_\xi =0 \ \mbox{in} \ \Sigma_L, \frac{\partial G_L}{\partial \nu}=0 \ \ \mbox{on} \ \partial \Sigma_L.\]
We have to show that
\begin{equation}
\label{A4}
G_L (x, y) \leq \frac{C}{|x-y|^{N-2}}
\end{equation}
where $C$ is independent of $L \geq L_0$,
 But (\ref{A4}) follows from standard potential estimates. See \cite{RW}.
 Note also that $ |f(y)|\leq \frac{C}{ (1+|y|)^{N-1}} < \frac{C}{ (1+ |y|)^3}$ so the integral $\int_{\Sigma_\ep} \frac{|f(y)|}{|x-y|^{N-2}} \leq C \|f\|_{**}  \frac{1}{(1+|y|)^{N-2}}$.

\bigskip

\noindent
{\bf Appendix B: Asymptotic behavior when $L \to L^{*}$}
\bigskip

In this appendix, we study the asymptotic behavior of the least energy solution when $L \to L^{*}$.  Let $ L=L_{*}$ and $ w_0$ be the unique radial solution of (\ref{ground}). It is clear that when $L \to L_{*}$, $ u_L \to w_0$ uniformly. In the following, we shall derive the next two order terms in the expansion of $ u_L$.

First, we consider the following linear problem
\begin{equation}
\Delta \phi -L_{*}^2 \phi + p w_0^{p-1} \phi=0 \ \mbox{in} \ \Sigma, \phi \in H^1 (\Sigma), \ \frac{\partial \phi}{\partial \nu}=0 \ \mbox{on} \ \Sigma
\end{equation}
Then  by separation of variables we have
\begin{equation}
\phi = \sum_{j=1}^{N-1} \frac{\partial w_0}{\partial y_j} + c_N \Phi_0
\end{equation}
where
\begin{equation}
\Phi_0= \phi_0 (|y^{'}|) \cos (\pi y)
\end{equation}
and $\phi_0$ is the principal eigenfunction of $ \Delta_{y^{'}} -L_{*}^2 + p w_0^{p-1}$. (We choose $\Phi_0$ so that $ \int_\Sigma \Phi_0^2=1$.)

As a consequence of the Fredholm Alternative, there exists a unique solution to the following problem
\begin{equation}
\Delta \phi -L_{*}^2 \phi + p w_0^{p-1} \phi=f \ \mbox{in} \ \Sigma, \phi \in H^1 (\Sigma), \ \frac{\partial \phi}{\partial \nu}=0 \ \mbox{on} \ \Sigma
\end{equation}
with
\begin{equation}
\int_{\Sigma} \phi \frac{\partial w_0}{\partial y_j}=\int_{\Sigma} \phi \Phi_0= \int_{\Sigma} f \frac{\partial w_0}{\partial y_j}=\int_{\Sigma} f \Phi_0=0
\end{equation}
such that
\begin{equation}
\| \phi \|_{H^2 (\Sigma)} \leq C \| f \|_{L^2 (\Sigma)}
\end{equation}

Let us denote $L_0:= \Delta -L_{*}^2  + p w_0^{p-1}$ and $ {\mathcal K}_0=\mbox{span} \ \{ \frac{\partial w_0}{\partial y_j}, \Phi_0 \}$. Now set
\begin{equation}
\delta = \int_{\Sigma} (u_L-w_0) \Phi_0
\end{equation}

We now claim that
\begin{equation}
\label{ul1}
u_L= w_0+ \delta \Phi_0 + \delta^2 \Phi_1 + \delta^3 \Phi_2
\end{equation}
where $\|\Phi_2\|_{H^2} \leq C$ and $\Phi_1$ satisfies
\begin{equation}
\label{A100}
L_0 \Phi_1+ \frac{L_{*}^2-L^2}{\delta^2} w_0 + \frac{p(p-1)}{2} w^{p-2} \Phi_0^2=0 \ \mbox{in} \ \Sigma, \Phi_1 \perp {\mathcal K}_0, \ \frac{\partial \Phi_1 }{\partial \nu}=0 \ \mbox{on} \ \Sigma.
\end{equation}

We prove it by expansion. Let $ u_L= w_0 + \delta \Phi_0 + \Phi_{1, L}$. Then $ \Phi_{1,L} \perp {\mathcal K}_0$ and $ \Phi_{1,L}= o(1)$. Since
$ \Delta (w_0+ \delta \Phi_0)-L^2 (w_0 +\delta \Phi_0)+ (w_0+\delta \Phi_0)^p=
  (L_{*}^2-L^2) w_0 + \delta^2 \frac{p(p-1)}{2} w^{p-2} \Phi_0^2 + {\mathcal O} (|\delta|^3 + |L_{*}^2-L^2| |\delta|)
$. We conclude that $ \Phi_{1,L}=   {\mathcal O} (|L_{*}^2-L^2|+ \delta^2)$. Now we let  $ \Phi_1$ be as defined in (\ref{A100}). Decomposing $u_l$ as in (\ref{ul1}), we see that $ \Phi_2$ satisfies
\[
L_0 [\Phi_2] + \frac{L_{*}^2-L^2}{\delta^2} \Phi_0 + p(p-1) w^{p-2} \Phi_0 \Phi_1 + \frac{p(p-1)(p-2)}{6} w^{p-3} \Phi_0^3 + O(\delta)=0
\]
and hence we have that
\[
 \frac{L_{*}^2-L^2}{\delta^2} \int_{\Sigma} \Phi_0^2 + p(p-1) w^{p-2} \int_{\Sigma} \Phi_0^2 \Phi_1 + \frac{p(p-1)(p-2)}{6}  \int_{\Sigma} w^{p-3} \Phi_0^4=0
\]
which gives the desired precise formula for $\delta$ in terms of $L_{*}^2-L^2$.

\end{document}